 \documentclass[final,3p,times]{elsarticle}
\usepackage{amsmath,amssymb}
\usepackage{natbib}
\usepackage{subfigure}                    %! used only when dealing with subfigures
\usepackage{multirow}
\journal{Non linear Analysis}

\newtheorem{theorem}{Theorem}
\newtheorem{proof}{Proof}

\newtheorem{definition}{Definition}
\begin{document}

\begin{frontmatter}



%\title{Sixth-order explicit one-step methods for stiff ODEs via hybrid deferred correction involving explicit second and forth order Runge-Kutta methods: Application to reaction-diffusion equations}

\title{Sixth-order explicit one-step methods for stiff ODEs via hybrid deferred correction involving RK2 and RK4: Application to reaction-diffusion equations}

%\author{Saint Cyr Elvi Rodrigue Koyaguerebo-Imé}
\author{Saint Cyr Elvi Rodrigue Koyaguerebo-Imé\corref{cor1}\fnref{label2}}
\ead{skoyague@univ-bangui.org}

\address{Université de Bangui, Département de Mathématiques et Informatique,  BP 1450, Bangui, Central African Republic.}

\begin{abstract}
In this paper, the fourth-order explicit Runge-Kutta method (RK4) is used to make a Deferred Correction (DC) on the explicit midpoint rule, resulting in an explicit one-step method of order six of accuracy, denoted DC6RK2/4. Convergence and order of accuracy of DC6RK2/4 are proven through a deferred correction condition satisfied by the RK4. The region of absolute stability of this method contains that of a RK6 and is tangent to the region [-5.626,0[x[-4.730,4.730] of the complex plane, containing a significant part of the imaginary axis. Numerical experiments with standard test problems for stiff systems of ODEs show that DC6RK2/4 performs well on problems regarding strong non-linearity and long-term integration, and this method does not require extremely small time steps for accurate numerical solutions of stiff problems. Moreover, this method is better than standard implicit methods like the Backward Differentiation Formulae and the DC methods for the implicit midpoint rule on stiff problems for which Jacobian matrices along the solution curve have complex eigenvalues where imaginary parts have larger magnitudes than real parts. An application of DC6RK2/4 to a class of test problems for reaction-diffusion equations in one dimensional is also carried out.
\end{abstract}


\end{frontmatter}




%% \bibliography{jabref3}


\section*{Introduction}

The aim of this paper is to study a correction of the explicit midpoint rule by the standard explicit fourth-order Runge-Kutta method (RK4) in order to derive an explicit one-step method that is expected to be of order six of accuracy and suitable for the numerical approximation of the initial value problem (IVP)
\begin{equation}
\label{a1}
\left\lbrace 
\begin{array}{cccc}
\displaystyle \frac{d u}{dt}&=&F(t,u),&~~ t\in [0,T],\\
u(0)&=&u_0,&
\end{array}
\right.
\end{equation}
where the unknown $u$ is from $[0,T]$ into a Banach space $X$, $u_0$ is a given data, and $F$ is a sufficiently differentiable function such that $u$ exists, is unique, and is sufficiently differentiable. This abstract Cauchy problem may result from the method of lines applied to Partial Differential Equations (PDEs).

 In fact, time-stepping methods constructed via deferred correction (DC) method, which consists of improving the order of accuracy of numerical methods of lower order, are widely investigated (see, e.g., \cite{gustafsson2001deferred,schild1990gaussian,auzinger2016encyclopedia,kushnir2012highly,hansen2011order,dutt2000spectral,
 daniel1967interated,IntegralDC2010,koyaguerebo2022arbitrary,bourgault2022variable}). In the case of time-stepping methods addressing stiff problems (problems extremely hard to
solve by standard explicit methods \cite{spijker1996stiffness}), a version of the DC method, which was introduced by Gustafsson and Kress in the reference \cite{gustafsson2001deferred}, is applied to the implicit midpoint rule, leading to an $A$-stable and $B$-convergent implicit scheme of order $2p+2$ of accuracy at the stage $p=0,1,\cdots$ of the correction (see \cite{koyaguerebo2022arbitrary}). The order of accuracy is guaranteed by a so-called deferred correction condition (DCC). Numerical experiments applied to different categories of Cauchy problems demonstrate the effectiveness of the DC method, and a complete analysis of the DC method, applied to reaction-diffusion equations, leads to an arbitrary high order and unconditionally stable method (see \cite{koyaguerebo2023arbitrary}). For non-stiff and moderately stiff problems, a number of schemes are constructed by applying iterative deferred correction (IDC) to lower-order methods such as the first-order Euler and Backward rules, the trapezoidal rule, the Taylor-2 and Adams-Bashforth of order 2 (see, e.g., \cite{auzinger2016encyclopedia,daniel1967interated,dutt2000spectral,hansen2011order,kushnir2012highly}) and to Runge-Kutta methods of order 2, 3, and 4 (see \cite{IntegralDC2010}). By IDC, the global error equation for an approximate solution of an IVP, obtained from a time-stepping method of order $p$ of accuracy, is transformed into a new IVP involving integral functions. The new IVP is solved by a time-stepping method of order $q$ of accuracy (a priori by the latter time-stepping method), using quadrature rules for the approximation of the integral functions, leading to a discrete approximation of the global error. The discrete approximation of the global error is locally added to the approximate solution of the basic IVP in order to successively improve its order of accuracy by $q$. In practice, IDC methods improve the region of absolute stability of schemes; however, resulting DC schemes are prone to high computational cost, order reduction, inaccuracy for large time steps and/or no asymptotic improvement for high levels of correction (see \cite{koyaguerebo2022arbitrary,koyaguerebo2023arbitrary} and references therein). Authors in the reference \cite{micalizzi2023new} have introduced interpolation processes in order to reduce the computational cost associated to IDC methods, but this strategy yields a poor improvement of the region of absolute stability. Besides IDC methods, there are the so-called explicit stabilized Runge-Kutta (ESRK) methods which are explicit one-step methods with extended stability domains along the negative real axis (see, e.g., \cite{assyr2010explicit, verwer1989convergence, assyr2002forth, skvortsov2011explicit, pintoa2023generalized}, and \cite{andrew2024family} for a version of two-step ESRK). These are a class of Runge-Kutta methods that use an arbitrary large internal stages by step via a methodical use of Chebychev polynomial or an extrapolation process, yielding explicit one-step methods with stability regions containing up to a segment $\left[ -\beta s^2,0\right] $ of the real axis (see, e.g., \cite{assyr2010explicit, verwer1989convergence}), where $s$ denotes the number of internal stages and $\beta$ is a positive real. Because of this property, the ESRK methods are considered to be appropriate for stiff ODEs resulting from space discretizations of parabolic or hyperbolic-parabolic equations (see \cite{verwer1989convergence, andrew2024family} and references therein). In practice, extended real negative stability interval of the ESRK methods is a good property for these methods to be efficient time-stepping methods for hyperbolic equations where exact solution has a lack of regularity since, in this case, lower order discretizations are needed in space and these require extremely small time step for standard explicit time-stepping methods (due to the constraint resulting from the CFL condition) or solutions of larger nonlinear equations by implicit methods. However, for any other parabolic or hyperbolic equations where exact solutions are sufficiently regular, a suitable discretization in space with high order methods is an ideal approach to avoid the need of large extended real stability interval. Nevertheless, when smaller tolerance is required or the right-hand side function $F$ in (\ref{a1}) is one-sided Lipschitz with a large Lipschitz constant or the magnitude of its derivatives is large, ESRK methods need to use small step sizes and then may become inefficient because of their larger internal stages that use too many function evaluations per step. We refer to \cite{skvortsov2011explicit} and references therein for inefficiency of high order ESRK methods regarding small tolerance and order reduction.

In general, stiff problems with Jacobian matrix having eigenvalues all negative along the solution curve behave as the ODE (\ref{a1}) where the function $F$ is one-sided Lipschitz with a Lipschitz constant $\mu \leq 0$. For such problems, local errors at a step $n$ with an implicit time-stepping method possesses, for example, a factor $(1-\mu k)^{-n}$ for the implicit Euler and $\left(\frac{2+\mu k}{2-\mu k}\right)^n$ for the implicit midpoint rule and its corrections. These factors allow to accurate approximations, even for time steps $k\geq 1$. However, when an explicit method with stability region containing $\mu$ is applied to such problems, the method is stable but not necessarily accurate for large time steps, requiring high order and smaller error constant. In the case where the Jacobian matrix has complex eigenvalues with dominant imaginary parts (this is the case for the problem B5 modified and the van der Pol oscillator in \cite{koyaguerebo2022arbitrary}), the extended stability domains along the negative real axis as $A$-stability are not sufficient to guarantee both stability and accuracy. We refer to \cite[Remark 4.2 \& Section 6.7]{koyaguerebo2022arbitrary} and \cite{stewart1990avoiding} for more details.

Our hybrid deferred correction scheme constitutes the first extension of the DC method introduced by Gustafsson and Kress in \cite{gustafsson2001deferred} and improved in \cite{koyaguerebo2022arbitrary, koyaguerebo2023arbitrary} to explicit time-stepping methods. The construction uses an approach of the deferred correction method that the current author and Yves Bourgault have developed in \cite{koyaguerebo2022arbitrary} for the computation of starting values (see \cite[Remark 2.1]{koyaguerebo2022arbitrary}). Applied to the explicit midpoint rule, the resulting scheme, for the numerical approximation of IVP (\ref{a1}), takes the simple form of the explicit midpoint scheme, but it contains two correction constants coming from finite difference approximations of local truncation errors of the explicit midpoint rule by approximate solutions of the IVP computed with RK4. Since this deferred correction method involves two different basic time-stepping methods, we call them hybrid deferred correction (an example of hybrid time-stepping method is presented in the reference \cite{ebadi2021new}). The hybrid approach is methodically chosen for two main reasons: (i) moderate the computational cost since the hybrid method of order six requires only twenty-one function evaluations per step while two successive corrections on RK2, in order to obtain an explicit DC method of order 6, require forty-two function evaluations per step; (ii) take advantage of the relatively larger region of absolute stability of $RK4$ and improve it.
As in \cite{koyaguerebo2022arbitrary}, convergence and order of accuracy of the hybrid DC schemes are guaranteed by a DCC that we define on RK4. Contrarily to the case of IDC methods, to apply the DC method developed in \cite{gustafsson2001deferred,koyaguerebo2022arbitrary} to any time-stepping method, knowledge of exact local truncation errors for this time-stepping method, written in terms of Taylor expansion up to an expected order of derivative, is required. This requirement disadvantages the application of this version of the DC method to correct RK4 to a method of order eight of accuracy. Numerical experiments with the hybrid deferred correction scheme for the explicit midpoint rule are performed on various IVPs regarding stiffness and a class of test problems related to one dimensional reaction-diffusion equations in order to assess its effectiveness. 
 
The paper is organized as follows: In Section \ref{sec:main}, we recall some basic results from finite difference and present the construction of the hybrid DC schemes. In Section \ref{sec:3}, we introduce a DCC for RK4 and prove a sufficient condition for RK4 to satisfy the DCC. In Section \ref{sec:4}, we present results about the convergence and order of accuracy of the hybrid DC methods. Section \ref{sec:5} deals with the absolute stability of the methods. Section \ref{sec:6} presents numerical experiments. Section \ref{sec:7} carries out applications of the new DC method to reaction-diffusion equation. Finally, in Section \ref{sec:8}, we comment on numerical results.

\section{Finite difference tools and construction of the hybrid deferred correction scheme}
\label{sec:main}
We suppose that $\displaystyle F\in C^{6}\left( [0,T]\times X,X \right)$ and such that the IVP (\ref{a1}) has a unique solution $\displaystyle u\in C^{7}\left( [0,T],X \right)$. We simply denote by $\|\cdot\|$, the norm in the Banach space $X$. For a time step $k>0$, we denote $t_n=nk$ and $t_{n+\lambda}=(n+\lambda)k$, for each integer $n$ and a positive real $\lambda$. This implies that $t_0=0$.  We consider time steps $k$ such that $0=t_0<t_1<\cdots< t_N=T$ is a partition of $[0,T]$, for a positive integer $N$. We have the following finite difference formulae about the centered and forward difference operators $D$ and $D_+$, respectively, associated to the step $k$ and applied to $u$:
$$
Du(t_{n+1/2})=\frac{u(t_{n+1})-u(t_n)}{k},$$
$$D_+u(t_{n})=\frac{u(t_{n+1})-u(t_n)}{k}.$$ 
Composites of $D_+$ are defined by

\begin{equation}
\label{bb2}
D_+^mu(t_n)=D_+\left(D_+^{m-1}u(t_n) \right) =k^{-m}\sum_{i=0}^{m}(-1)^i{{m}\choose {i}}u(t_{n+m-i}),
\end{equation}
for any integer $m\geq 1$. We have the estimate
\begin{equation}
\label{bb3}\left \Vert D_+^{m}u(t_n) \right \Vert\leq \max_{0\leq t\leq  T}\left \Vert\frac{d^{m}u}{dt^{m}}(t)\right \Vert,
\end{equation}
provided $[t_{n},t_{n+m}] \subset [0,T]$ and $m\leq p+1$ (see, e.g.,\cite[Theorem 4]{koyaguerebo2022finite} or \cite[p.249]{isaacson1966analysis} ). If $\left\lbrace u^n \right\rbrace_n$ is a sequence of approximations of $u$ at the discrete points $t_n$, the finite difference operators apply to $\left\lbrace u^n \right\rbrace _n$, and we define
$$
Du^{n+1/2}=D_+u^{n}=D_-u^{n+1}=\frac{u^{n+1}-u^n}{k}.$$

 Let us consider the following Taylor expansions of $u$ about the points $t_{n+1/2}$ and $t_n$, respectively:
\begin{equation*}
u(t_{n+1})=u(t_n)+ku'(t_{n+\frac{1}{2}})+\sum_{i=1}^3\frac{2k^{2i+1}}{2^{2i+1}\times (2i+1)!}u^{(2i+1)}(t_{n+\frac{1}{2}})+O(k^9),
\end{equation*}
and
\begin{equation*}
u(t_{n+\frac{1}{2}}) =u(t_{n})+\frac{k}{2}u'(t_n)+\sum_{i=2}^6\frac{k^i}{2^i\times i!}u^{(i)}(t_n)+O(k^{7}).
\end{equation*} 
We can rewrite these two expansions as follows: 
\begin{equation}
\label{b1}
u(t_{n+1})=u(t_n)+ku'(t_{n+\frac{1}{2}})+\sum_{i=1}^3\frac{2\times 5^{2i+1} h^{2i+1}}{2^{2i+1}\times (2i+1)!}u^{(2i+1)}(t_{n+\frac{1}{2}})+O(k^9)
\end{equation}
and
\begin{equation}
\label{b2}
u(t_{n+\frac{1}{2}}) =u(t_{n})+\frac{k}{2}u'(t_n)+\sum_{i=2}^6\frac{5^i\times h^i}{2^i\times i!}u^{(i)}(t_n)+O(k^{7}),
\end{equation}  
where $h=k/5$. We substitute the following identity in (\ref{b1})
$$h^3u^{(3)}\left( t_{n+\frac{1}{2}}\right)= h^3(D_+D_-)Du(t_{n+\frac{1}{2}})-\frac{480}{2^5\times 5!}h^5u^{(5)}\left( t_{n+\frac{1}{2}}\right)-\frac{4368}{2^7\times 7!}h^7u^{(7)}+O(h^9),$$ where
$$h^3(D_+D_-)Du(t_{n+1/2})=u(t_n+4h)-3u(t_n+3h)+3u(t_n+2h)-u(t_n+h).$$
This yields
\begin{equation*}
\begin{aligned}
u(t_{n+1})=u(t_n)+&ku'(t_{n+\frac{1}{2}})+\frac{125}{24}h^3(D_+D_-)Du(t_{n+\frac{1}{2}})+\frac{125}{128}h^5u^{(5)}\left( t_{n+\frac{1}{2}}\right)+\frac{133500}{2^7\times 7!}h^7u^{(7)}(t_{n+\frac{1}{2}})+\cdots.
\end{aligned}
\end{equation*}
Substituting the identity
$$
\begin{aligned}
h^5(D_+D_-)^2Du(t_{n+\frac{1}{2}})&=\sum_{j=0}^{5}(-1)^j\binom{5}{j}u( t_n+(5-j)h)=h^5u^{(5)}(t_{n+\frac{1}{2}})+\frac{134400}{2^7\times 7!}h^7u^{(7)}(t_{n+\frac{1}{2}})+\cdots
\end{aligned}
$$
in the last one, we deduce the discrete expansion
\begin{equation}
\label{b3}
\begin{aligned}
u'(t_{n+\frac{1}{2}})=\frac{u(t_{n+1})-u(t_n)}{k}-\frac{125}{24k}h^3(D_+D_-)Du(t_{n+\frac{1}{2}})-\frac{125}{128k}h^5(D_+D_-)^2Du(t_{n+\frac{1}{2}})-\frac{2250}{10^7\times 7!}k^7u^{(7)}(t_{n+\frac{1}{2}})+O(k^9).\end{aligned}
\end{equation}
Similarly, substituting successively the identities
\begin{equation*}
\begin{aligned}
h^mu^{(m)}(t_n)=h^mD_+^mu(t_n)-\sum_{i=m+1}^{\infty} \frac{k^{i}}{i!}u^{(i)}(t_{n})\sum_{j=0}^{m}(-1)^j\binom{m}{j}(m-j)^{i},
\end{aligned}
\end{equation*}
$m=2,3,4,5$, in (\ref{b2}), we obtain the discrete expansion
\begin{equation}
\label{b4}
\begin{aligned}
u&(t_{n+1/2}) =u(t_{n})+\frac{k}{2}u'(t_n)+\frac{25}{8}h^2D_+^2u(t_n)-\frac{25}{48}h^3D_+^3u(t_n)\\&+\frac{225}{2^4\times 4!}h^4D_+^4u(t_n)-\frac{1875}{2^5\times 5!}h^5D_+^5u(t_n)+\frac{18975}{10^6\times6!}k^6u^{(6)}(t_n)+O(k^7).
\end{aligned}
\end{equation}
Therefore, evaluating the first equation in (\ref{a1}) at $t_{n+1/2}$, we deduce from (\ref{b3}) and (\ref{b4}) that

\begin{equation}
\label{b5}
\left\lbrace \begin{aligned}
 \displaystyle \frac{u(t_{n+1})-u(t_n)}{k}-\frac{1}{k} E_1( t_{n+\frac{1}{2}})& = F\left( t_{n+\frac{1}{2}},u(t_{n})+\frac{k}{2}F(t_n,u(t_n))+E_2( t_{n+\frac{1}{2}})\right),\\
 u(0)&=u_0,
\end{aligned}\right. 
\end{equation}
where
$$
\begin{aligned}
E_1\left( t_{n+\frac{1}{2}}\right)&=\frac{125}{24}h^3(D_+D_-)Du(t_{n+\frac{1}{2}})+\frac{125}{128}h^5(D_+D_-)^2Du(t_{n+\frac{1}{2}})+O(k^7)\\&=\frac{125}{384}\left[-3u(t_n)-u(t_{n+\frac{1}{5}}) +18u(t_{n+\frac{2}{5}})-18u(t_{n+\frac{3}{5}})+u(t_{n+\frac{4}{5}})+3u(t_{n+1})\right] +O(k^7),
\end{aligned}
$$
and
$$
\begin{aligned}
&E_2( t_{n+\frac{1}{2}})=\frac{25}{8}h^2D_+^2u(t_n)-\frac{25}{48}h^3D_+^3u(t_n)+\frac{75}{128}h^4D_+^4u(t_n)-\frac{125}{256}h^5D_+^5u(t_n)+O(k^6).\end{aligned}$$
%(25/768)*[145;-387;402;-238;93;-15]
It follows the scheme
\begin{equation}
\label{b6}
\left\lbrace 
\begin{aligned}
\displaystyle 
u^{6,0}&=u_0,\\
u^{4,n,0}&=u^{6,n},\\
a_n=&\frac{125}{384}\left(-3u^{4,n,0}-u^{4,n,1}+18u^{4,n,2}-18u^{4,n,3}+u^{4,n,4}+3u^{4,n,5}\right), \\
b_n=&\frac{25}{768}(145u^{4,n,0}-387u^{4,n,1}+402u^{4,n,2}-238u^{4,n,3}+93u^{4,n,4}-15u^{4,n,5}),\\
u^{6,n+1}&=u^{6,n}+a_n+kF\left( t_{n+\frac{1}{2}},u^{6,n}+\frac{k}{2}F(t_n,u^{6,n})+b_n\right).
\end{aligned}
\right.
\end{equation}
The sequence $\left\lbrace u^{6,n}\right\rbrace _{n=0}^N$ is expected to be an approximation of the exact solution $u$ of the IVP (\ref{a1}) with order 6 of accuracy at the corresponding points $t_n$, and quantities $u^{4,n,i}$ are approximations of $u(t_{n+i/5})$, $i=0,1,2,\cdots,5$, computed successively with the scheme RK4 at the points $t_{n+i/5}$ for each fixed $n$. Quantities $a_n$ and $b_n$ are approximations of the error terms $E_1\left( t_{n+1/2}\right)$ and $E_2\left( t_{n+1/2}\right)$, respectively, by the approximate solutions $u^{4,n,i}$. We call the scheme (\ref{b6}) hybrid deferred correction involving RK2 and RK4, denoted $DC6RK2/4$.

\section{Deferred correction condition for the fourth-order explicit Runge-Kutta method}
\label{sec:3}
In this section, we give a deferred correction condition (DCC) for RK4 and establish a DCC property of $RK4$, which is necessary for the proof of convergence and order of accuracy of the scheme $DC6RK2/4$ in the next section. Hereafter, the letter $C$ will denote any constant independent of $k$ and that can be calculated explicitly in terms of known quantities. The exact value of $C$ may change from a line to another. We have the following definition:

\begin{definition}[Deferred Correction Condition for RK4]
\label{def:1} Let $u$ be the exact solution of the Cauchy problem (\ref{a1}). A sequence $\left\lbrace u^{4,n}\right\rbrace _{n=0}^N$ of approximations of $u$, at the discrete points $0=t_0<\cdots <t_N=T$, is said to satisfy Deferred Correction Condition $(DCC)$ for the fourth-order explicit Runge-Kutta method (RK4) if $\left\lbrace u^{4,n}\right\rbrace _{n=0}^N$ approximates $u(t_n)$ with order 4 of accuracy, and we have
\begin{equation}
\label{b7} 
 \Vert D_+^i\left( u^{4,n}-u(t_{n})\right) \Vert \leq C k^{4},\quad 
\end{equation}
for each integer $i=0,1,2,\cdots$, and $n=0,1,2,...,N-i$, provided $u$ is $C^{5+i}$, where $C$ is a constant independent of $k$.
\end{definition}

The following theorem gives a sufficient condition for RK4 to satisfy DCC.
\begin{theorem}
\label{thm:2}Suppose that $\displaystyle F\in C^{p}\left( [0,T]\times X,X \right)$, $p\geq 6$, and such that problem (\ref{a1}) has a unique solution $\displaystyle u\in C^{p+1}\left( [0,T],X \right)$. 
Let $\left\lbrace u^{4,n}\right\rbrace _{n=0}^N$ be the sequence of approximate solutions of (\ref{a1}) by the fourth order explicit Runge-Kutta method on the partition $0=t_0<\cdots <t_N=T$, $t_n=nk$, $k=T/N$, of $[0,T]$. Then $\left\lbrace u^{4,n}\right\rbrace _{n=0}^N$ satisfies the deferred correction condition for the fourth-order explicit Runge-Kutta method: 
\begin{equation}
\label{b8} 
 \Vert D_+^m\left( u^{4,n}-u(t_{n})\right) \Vert \leq C k^{4},\;\; m=0,1,2,\cdots,p-4, n=0,1,2,...,N-m,
\end{equation}
where $C$ is a constant independent of $k$.
\end{theorem}

\begin{proof}For the sake of simplification, we suppose that $F=F(x)$. The general case can be deduced from elementary (but tedious) calculations. 

Since $\left\lbrace u^{4,n}\right\rbrace _{n=0}^N$ is a sequence of approximate solutions of (\ref{a1}) by RK4 on the uniform partition $0=t_0<\cdots <t_N=T$, $t_n=nk$, $k=T/N$, of $[0,T]$, we have
\begin{equation}
\label{b9}
\frac{u^{4,n+1}-u^{4,n}}{k}=\frac{1}{6}\left(K_1^n+2K_2^n+2K_3^n+K_4^n \right),
\end{equation}
where
$$
\begin{aligned}
K_1^n=F\left( u^{4,n}\right) ;\quad K_2^n=F\left( u^{4,n}+\frac{k}{2}K_1^n\right);\quad K_3^n=F\left( u^{4,n}+\frac{k}{2}K_2^n\right); \;K_4^n=F\left(  u^{4,n}+k K_3^n\right) .
\end{aligned}
$$
Otherwise, we can write
\begin{equation}
\label{b10}
\frac{u(t_{n+1})-u(t_n)}{k}=\frac{1}{6}\left(\kappa_1(t_n)+2\kappa_2(t_n)+2\kappa_3(t_n)+\kappa_4(t_n) \right)+k^4\tau(t_n),
\end{equation}
where $\tau$ is $ C^{p-4}\left( [0,T],X \right)$,
$$\kappa_1(t_n)=F\left( u(t_n)\right) ;\kappa_2(t_n)=F\left( u(t_n)+\frac{k}{2}\kappa_1(t_n)\right);$$
$$ \kappa_3(t_n)=F\left( u(t_n)+\frac{k}{2}\kappa_2(t_n)\right); \kappa_4(t_n)=F\left( u(t_n)+k\kappa_3(t_n)\right) .$$
Hypotheses of Theorem \ref{thm:2} implies that $\left\lbrace u^{4,n}\right\rbrace _{n=0}^N$ approximates $u$ with order four of accuracy:
\begin{equation}
\label{b11}\Vert u(t_n)-u^{4,n} \Vert \leq Ck^4, \mbox{ for each } n=0,1,2,\cdots, N,
\end{equation}where $C$ is a constant depending only on $T$, $F$ and derivatives of $u$ up to order 5.
To establish (\ref{b8}), we proceed by induction on the integer $m=0,1,2,..., p$.
\newline

\noindent
(i)\;\; Inequality (\ref{b8}) for $m=0$. This is done by inequality (\ref{b11}).
\newline

\noindent
(ii)\;\; Here we are going to prove that the inequality (\ref{b8}) remains true for $m+1$, assuming that it holds for an arbitrary integer $m$ such that $0\leq m\leq p-1$.

\noindent
We combine (\ref{b9}) and (\ref{b10}) and deduce the identity
\begin{equation}
\label{a44}
\begin{aligned}
D_+\Theta^{4,n}=&\frac{1}{6}\left[(K_1^n-\kappa_1(t_n))+2(K_2^n-\kappa_2(t_n))+2(K_3^n-\kappa_3(t_n))+(K_4^n-\kappa_4(t_n) )\right]-h^4\tau(t_n),
\end{aligned}
\end{equation}where 
 $$\Theta^{4,n}=u^{4,n}-u(t_n).$$
Let
$$h_1(t_n)=K_1^n-\kappa_1(t_n)=\int_0^1dF(A_1^n)\left(\Theta^{4,n} \right)\,d\tau_1, $$
where
$$A_1^n=u(t_n)+\tau_1\Theta^{4,n}.$$
Setting 
$$d\tau^i =d\tau_1\cdots d\tau_i,$$and
\begin{equation}
\label{a13}
A_{i+1}^{n}=A^n_{i}+\tau_{i+1}(A_{i}^{n+1}-A_{i}^{n})=A_1^n+\sum_{l=1}^i \sum_{2\leq i_1<\cdots <i_l \leq i+1}\tau_{i_1}\cdots\tau_{i_l}k^lD_+^lA_1^n,
\end{equation}
for each positive integer $i$, we have
$$D_+h_1(t_{n})=\int_0^1\int_0^1d^2F(A_2^{n})(D_+A_1^n,{\Theta}^{4,n+1})d\tau^2 +\int_0^1dF\left(A_1^n\right)( D_+{\Theta}^{4,n}) d\tau_1,$$
and
\begin{equation*}
\begin{aligned}
&D^2_+h_1(t_{n})=\int_0^1\int_0^1\int_0^1d^3F(A_3^n)\left(D_+A_2^n, D_+A_1^{n+1},{\Theta}^{4,n+2}\right) d\tau^3+\int_0^1\int_0^1d^2F(A_2^n)\left(D_+^2A_1^n,{\Theta}^{4,n+2}\right) d\tau^2\\&+2\int_0^1\int_0^1d^2F(A_2^n)\left(D_+A_1^n,D_+{\Theta}^{4,n+1}\right) d\tau^2+\int_0^1dF(A_1^n)(D^2_+{\Theta}^{4,n})d\tau_1.
\end{aligned}
\end{equation*}
More generally, we have
\begin{equation}
\label{a14}
\begin{split}
D_+^qh_1(t_n)=\sum_{i=1}^{q+1}\sum_{|\alpha_i|=q} a_{\alpha_i}L^{n,q}_{\alpha_i} ,\mbox{ for } q=1,2,...,p, \mbox{ and }n\leq N-q,
\end{split}
\end{equation}
where $\alpha_i=(\alpha_i^1,\alpha_i^2,\cdots,\alpha_i^i) \in \left\lbrace 0,1,2,\cdots, i-1\right\rbrace^i $, $\vert \alpha_i\vert =\alpha_i^1+\alpha_i^2+\cdots+\alpha_i^i$, the $a_{\alpha_i}$ are constants, and

$$L^{n,q}_{\alpha_i}=\int_{[0,1]^i}d^iF(A_i^{n})(D_+^{\alpha_i^i}A_{i-1}^{n+r_i^i},D_+^{\alpha_i^{i-1}}A_{i-2}^{n+r_i^{i-1}},\cdots,D_+^{\alpha_i^2}A_{1}^{n+r_i^2},D_+^{\alpha_i^1}{\Theta}^{4,n+r_i^1})d\tau^i ,$$
with $l+r_i^l=i$ and $\alpha_i^l+r_i^l\leq q$, for $l=1,2,\cdots,i$. From the definition of $A_1^n$ and the identity (\ref{a13}), we have 
$$A_i^{n+1}=u(t_{n})+O(k^4), \mbox{ for } i=1,2,\cdots,p, $$ and we deduce from the regularity of the function $F$ that
\begin{equation}
\label{a50}\left \Vert d^i F\left( A_{i}^{n} \right) \right \Vert \leq C_i,\;\mbox{ for } i=1,2,...,p, \mbox{ and }0\leq n\leq N-p,
\end{equation}
where $C_i$ is a constant depending only on $T$ and the derivatives of $F$ and $u$ up to order $i$ and $i+1$, respectively. Otherwise, from the induction hypothesis and the continuity of $u$, we have
\begin{equation}
\label{a51}\Vert D_+^{\alpha_i^{l}}A_{l-1}^{n+r_i^{l}}\Vert \leq C, \mbox{ for } 2\leq l \leq i-1\leq m,
\end{equation}
and 
\begin{equation}
\label{a52}
\Vert D_+^{\alpha_i^1}{\Theta}^{4,n+r_i^1}\Vert \leq Ck^4, \mbox{ for } i \leq m+1,
\end{equation}
 knowing that $\alpha_i^1\leq i-1$, where $C$ is a constant depending only on $m$, $T$ and the derivatives of $F$ and $u$ up to order $i$ and $i+1$, respectively. Each $L^{n,q}_{\alpha_i}$ being multilinear continuous function, we deduce from (\ref{a50})-(\ref{a52}) that
$$\Vert L^{n,q}_{\alpha_i}\Vert \leq Ck^4 ,\mbox{ for } 1\leq i\leq m.$$
It follows by the triangle inequality and identity (\ref{a14}) for $q=m$ that
\begin{equation*}
\Vert D_+^{m} h_1(t_{n}) \Vert \leq Ck^4,\quad n=0,1,2,\cdots, N-(m+1),
\end{equation*}
where $C$ is a constant depending only on $p$, $T$ and the derivatives of $F$ and $u$ up to order $m+1$ and $m+2$, respectively. The last inequality implies that,
\begin{equation}
\label{b18}\Vert D_+^{m} \left( K_1^n-\kappa_1(t_{n})\right)  \Vert \leq Ck^4,\quad n=0,1,2,\cdots, N-(m+1).
\end{equation} 
On the other hand, we can write
$$h_j(t_n)=K_j^n-\kappa_j(t_n)=\int_0^1dF(A_{1,j}^n)\left(A_{0,j}^{n} \right)\,d\tau_1, \quad j=2,3,4,$$
with $$A_{0,j}^{n}={\Theta}^{4,n}+a_jk\tau_1\left( K_{j-1}^n-\kappa_{j-1}(t_n)\right), \quad a_4=2a_2=2a_3=1,$$
and
$$A_{1,j}^{n}=K_{j-1}^n+\tau_1\left( K_{j-1}^n-\kappa_{j-1}(t_n)\right).$$
Therefore, proceeding as in (\ref{b18}), replacing $A_j^n$ by $A_{1,j}^n$ and ${\Theta}^{4,n}$ by $A_{0,j}^n$, we successively deduce that
\begin{equation}
\label{b19}\Vert D_+^{m} \left( K_j^n-\kappa_j(t_{n})\right)  \Vert \leq Ck^4,\quad n=0,1,2,\cdots, N-(m+1),\; j=2,3,4.
\end{equation} 
Therefore, applying $D_+^m$ to (\ref{a44}), we deduce from the triangle inequality, inequalities (\ref{b18})-(\ref{b19}) and the estimate (\ref{bb3}) applied to the function $\tau$ that
$$\Vert D_+^{m+1}\Theta^{2,n}\Vert \leq Ck^4,$$
for a constant $C$ independent of $k$. Finally, by induction, Theorem \ref{thm:2} is proven.
\end{proof}

\section{Convergence and order of accuracy}
\label{sec:4}
In this section, we propose to establish the convergence and order of accuracy of the hybrid DC scheme (\ref{b6}). The proof relies on the DCC related to the explicit $RK4$, given in Section \ref{sec:3}. We have the following result:

\begin{theorem}
\label{thm:3}Suppose that $\displaystyle F\in C^{6}\left( [0,T]\times X,X \right)$ and Lipschitz with respect to its second variable such that problem (\ref{a1}) has a unique solution $\displaystyle u\in C^{7}\left( [0,T],X \right)$.
More precisely, we suppose that there exists $\mu \geq 0$ such that
\begin{equation}
\label{a31a} \|F(t,x)-F(t,y)\|\leq \mu\|x-y\|,~~\forall (t,x,y)\in [0,T]\times X\times X.
\end{equation}
%
%\item[(ii)] $X$ is a Hilbert space with inner product $\left( .,. \right) $, and there exists $\mu \geq 0$ such that 
%\begin{equation}
%\label{a31} \left( F(t,x)-F(t,y),x-y\right)\leq \mu\|x-y\|^2, ~~\forall (t,x,y)\in [0,T]\times X\times X. 
%\end{equation}
%\end{enumerate}
Then the sequence $\left\lbrace u^{6,n}\right\rbrace _{n=0}^N$, given by the scheme (\ref{b6}), approximates $u(t_n)$ with order 6 of accuracy on the partition $0=t_0<\cdots <t_N=T$, $t_n=nk$, $k=T/N$, of $[0,T]$, provided $k$ is sufficiently small.
\end{theorem}

\begin{proof}
Combining the first identity in (\ref{b5}) and the last identity in (\ref{b6}), we obtain
\begin{equation}
\label{b24a}\begin{aligned}
&u^{6,n+1}-u(t_{n+1})=u^{6,n}-u(t_{n})+k F\left( t_{n+\frac{1}{2}},u^{6,n}+\frac{k}{2}F(t_n,u^{6,n})+b_n\right)\\&- kF\left( t_{n+\frac{1}{2}},u(t_n)+\frac{k}{2}F\left( t_n,u(t_n)\right) +E_2\left( t_{n+\frac{1}{2}}\right) \right) +a_n-E_1\left( t_{n+\frac{1}{2}}\right).
\end{aligned}
\end{equation}
 We are going to prove by induction on the integers $n=0,1,2,\cdots,N$ that
 \begin{equation}
 \label{b24}
 \Vert u^{6,n+1}-u(t_{n+1})\Vert \leq \left(1+\mu k+\frac{1}{2}\mu^2 k^2 \right)\Vert u^{6,n}-u(t_{n})\Vert+Ck^7,
 \end{equation}
where $C$ is a constant independent of $n$. 
 
 \vspace{0.5cm}
 \noindent
 (i)\; The case $n=0$.
 
 From the hypotheses of Theorem \ref{thm:3}, $\left\lbrace u^{4,0,i}\right\rbrace_i $ satisfies DCC for $RK4$. It follows from the triangle inequality that
 $$
\begin{aligned} 
 \Vert a_0-E_1(t_{\frac{1}{2}})\Vert &\leq \frac{125}{4}h^3\Vert D_+^3\Theta^{4,0,0}\Vert+\frac{125}{128}h^5\Vert D_+^5\Theta^{4,0,0}\Vert+Ch^7\leq C(1+C_{DCC}) h^7,
\end{aligned} 
 $$
 where $\Theta^{4,n,i}=u^{4,n,i}-u(t_n+ih)$, $n=0,1,2,\cdots,N$ and $i=0,1,\cdots,5$, $C_{DCC}$ is the constant due to the deferred correction condition in Theorem \ref{thm:2}, and the finite defference operator $D_+$ is applied to the third argument $i$ of  $\Theta^{4,n,i}$. Proceeding similarly, we obtain
 $$\Vert b_0-E_2(t_{\frac{1}{2}})\Vert \leq Ch^6,$$
Therefore, passing to the norm $\Vert \cdot\Vert$ on both sides of the last identity for $n=0$, we deduce from the last two inequalities and the Lipschitz condition (\ref{a31a}) that
 $$\Vert u^{6,1}-u(t_{1})\Vert \leq \left(1+\mu k+\frac{1}{2}\mu^2 k^2 \right)  \Vert u^{6,0}-u_0)\Vert+Ck^7\leq Ck^7,$$where $C$ is a constant independent of $k$.

 \vspace{0.5cm}
 \noindent
 (ii)\; We prove (\ref{b24}) for $n+1$, assuming that this inequality is true for an arbitrary nonegative integer $n$.
 
From the hypotheses of Theorem \ref{thm:3}, $\left\lbrace u^{4,n+1,i}\right\rbrace_i $ satisfies the DCC for $RK4$ since it is obtained from the explicit RK4 with the starting value $u^{4,n+1,0}=u^{6,n+1}$ which is of order six of accuracy, according to the induction hypothesis. It follows from the triangle inequality that
 $$
\begin{aligned} 
 \Vert a_{n+1}-E_1(t_{n+\frac{3}{2}})\Vert &\leq \frac{125}{4}h^3\Vert D_+^3\Theta^{4,n+1,0}\Vert+\frac{125}{128}h^5\Vert D_+^5\Theta^{4,n+1,0}\Vert+Ch^7\leq C(1+C_{DCC})h^7,
\end{aligned} 
 $$and 
 $$\Vert b_{n+1}-E_2(t_{n+\frac{3}{2}})\Vert \leq C(1+C_{DCC})h^6.$$
Therefore, passing to norm $\Vert \cdot\Vert$ on both sides of (\ref{b24a}), we deduce from the last two inequalities and the Lipschitz condition (\ref{a31a}) that
 $$\Vert u^{6,n+2}-u(t_{n+2})\Vert \leq \left(1+\mu k+\frac{1}{2}\mu^2 k^2 \right)  \Vert u^{6,n+1}-u(t_{n+1})\Vert+Ck^7,$$ where $C$ is a constant independent of $k$.

\noindent
 Form (i) and (ii), we deduce by induction that (\ref{b24}) is true for any integer $n\geq 0$. It follows from (\ref{b24}) that
 \begin{equation*} 
 \Vert u^{6,n}-u(t_{n})\Vert \leq \left(1+\mu k+\frac{1}{2}\mu^2 k^2 \right)^n\Vert u^{6,0}-u(t_{0})\Vert+nCk^7.
 \end{equation*}
Since $T=Nk>nk$ and $u^{6,0}=u(t_0)=u_0$, the last inequality yields
\begin{equation}
 \label{b26}
 \Vert u^{6,n}-u(t_{n})\Vert \leq Ck^6,
 \end{equation}
 for each integer $n=0,1,2,\cdots,N$, where $C$ is a constant depending only on $T$ and the derivatives of $u$ and $F$. Whence, $\left\lbrace u^{6,n}\right\rbrace _{n=0}^N$ approximates $u(t_n)$ with order 6 of accuracy.
\end{proof}

\section{Region of absolute stability for the hybrid DC scheme}
\label{sec:5}
\begin{figure}[!h]%% placement specifier
%% Use \includegraphics command to insert graphic files. Place graphics files in 
%% working directory.
\centering%% For centre alignment of image.
\includegraphics[height=2.77in, width=3.5in]{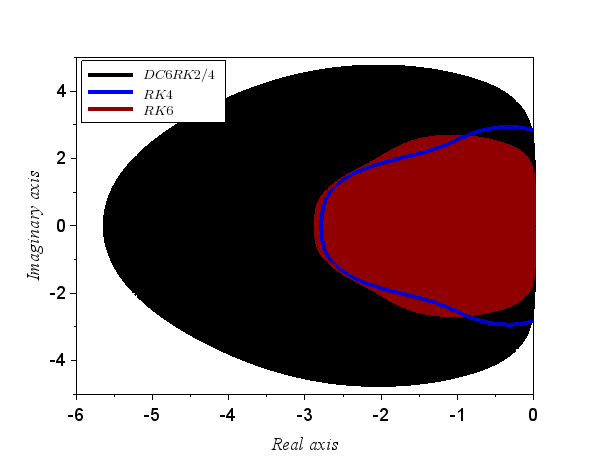}
%% Use \caption command for figure caption and label.
\caption{Regions of absolute stability for $DC6RK2/4$, $RK6$ and $RK4$.}
\label{Fig}
%% https://en.wikibooks.org/wiki/LaTeX/Importing_Graphics#Importing_external_graphics
\end{figure}

Applying the hybrid DC scheme (\ref{b6}) to the linear IVP 
\begin{equation}
\label{b27}
\left\lbrace 
\begin{array}{cccc}
\displaystyle \frac{d u}{dt}&=&\lambda u,&~~ t\in [0,T],\\
u(0)&=&1,&
\end{array}
\right.
\end{equation}
and setting $z=\lambda k$, we obtain
\begin{equation}
\left\lbrace 
\begin{aligned}
\displaystyle 
u^{6,0}&=1,\\
u^{4,n,i}&=q\left( \frac{z}{5}\right)^iu^{6,n},\quad i=0,1,\cdots, 5\\
a_n&=r(z)u^{6,n},\\
%a_n&=\frac{125}{384}\left(-3u^{4,n,0}-u^{4,n,1}+18u^{4,n,2}-18u^{4,n,3}+u^{4,n,4}+3u^{4,n,5}\right), \\
%b_n&=\frac{25}{768}(145u^{4,n,0}-387u^{4,n,1}+402u^{4,n,2}-238u^{4,n,3}+93u^{4,n,4}-15u^{4,n,5})\\
b_n&=s(z)u^{6,n},\\
u^{6,n+1}&=\left(1+z+\frac{z^2}{2}+ r(z)+zs(z)\right) u^{6,n},
\end{aligned}
\right.
\end{equation}
where $$q(z)=1+z+\frac{z^2}{2!}+\frac{z^3}{3!}+\frac{z^4}{4!}$$is the stability function of the explicit fourth-order Runge-Kutta method (see, e.g., \cite[page 202]{lambert1991numerical}),
$$r(z)=\frac{125}{384}\left(-3-q\left( \frac{z}{5}\right)+18q\left( \frac{z}{5}\right)^2-18q\left( \frac{z}{5}\right)^3+q\left( \frac{z}{5}\right)^4+3q\left( \frac{z}{5}\right)^5\right),$$
and
$$s(z)=\frac{25}{768}\left( 145-387q\left( \frac{z}{5}\right)+402q\left( \frac{z}{5}\right)^2-238q\left( \frac{z}{5}\right)^3+93q\left( \frac{z}{5}\right)^4-15q\left( \frac{z}{5}\right)^5\right).$$
It results the following stability function of the hybrid deferred correction scheme (\ref{b6}), which a polynomial of degree 21,
$$R(z)=1+z+\frac{z^2}{2}+ r(z)+zs(z).$$
The corresponding region of absolute stability together with that of $RK4$ and the sixth-order explicit Runge-Kutta method introduced by the author H. A. Luther in \cite{luther1968explicit} is presented in Figure \ref{Fig}.

\section{Numerical experiments}
\label{sec:6}
In this section, we evaluate the effectiveness of the scheme DC6RK2/4, implemented using the Scilab programming language. We consider six test problems related to ODEs, while an application to PDEs follows in the next section. The six test problems result from \cite{koyaguerebo2022arbitrary}, and they address $B$-convergence (problems \ref{b77b} and \ref{69}), long-term integration (problems \ref{66} and \ref{69}), and general stiffness.

For a comparison of accuracy of the scheme DC6RK2/4 with some existing time-stepping methods, we implement the explicit Runge-Kutta method of order 6 ($RK6$) introduced by the author H. A. Luther in \cite{luther1968explicit} and the standard $RK4$, and we consider $DC6$ (the deferred correction scheme of order 6 constructed on the implicit midpoint rule) and the Backward Differentiation Formula of order six ($BDF6$) which are evaluated on these test problems in \cite{koyaguerebo2022arbitrary}. In the case of problems without analytic solutions, we use the functions \texttt{stiff} (implementing BDF with adaptive steps) from Scilab to compute a reference solution and the starting values of the implemented BDF6. We recall that, applied to each of these six test problems, even with default tolerance, the solver \texttt{rkf} from Scilab reports that "it is likely that rkf45 is inefficient for solving this problem" or "integration was not completed because requested accuracy could not be achieved using smallest allowable stepsize".

For each of the problems, we give a table of absolute errors and orders of accuracy for pairs of two consecutive time steps. We use dashes in the tables to indicate absolute errors that are not a number or larger than $10^{16}$. Given an exact solution $u=(u_1,\cdots,u_d)~:~[0,T] \rightarrow \mathbb{R}^d$, $1\leq d\leq 6$, the absolute errors on the approximate solutions $\left\lbrace u^{n}_i \right\rbrace_{0\leq n\leq N, 1\leq i\leq d} $ is computed with the norm
$$\|\left\lbrace u^{n}_i\right\rbrace_n-u_i\| =\max _{0\leq n\leq N}|u^{n}_i-u_i(t_n)|,\quad 1\leq i\leq d.$$ 
For larger $N$, we extract solutions on at least $6\times 10^4$ discrete times evenly spread over the interval $[0,T]$.

\subsection{Bernoulli differential equation (see \cite{koyaguerebo2022arbitrary}), stiff with real negative eigenvalues}
\begin{equation}
\label{b77b}
u'(t)=F(t,u)=-0.1u(t)-1000u^{20}(t)~,~~~u(0)=1, ~~t\in [0,10].
\end{equation}
This problem was introduced by authors in \cite{koyaguerebo2022arbitrary} to address $B$-convergence. Table \ref{tab:1} gives the absolute error and the order of convergence for each pair of consecutive time steps. Dashes for RK6 indicate absolute errors superior or equal to $10^{16}$. The absolute errors for $BDF1$ and $BDF2$ when $k=4\times 10^{-3}$ are both $3.109 \times 10^{-2}$. 

From Table \ref{tab:1}, we observe that the implicit methods $DC6$ and $BDF6$ are more suitable for this problem in the case of large time steps. The method $DC6RK2/4$ achieves its proper order on this problem and produces accurate approximate solutions for time steps slightly large. However, $RK6$ requires smaller time steps for accurate approximations, but this method is more accurate than all the other methods in the region of asymptotic convergence $k\leq 10^{-4}$.
\begin{table}[!ht]
\caption{Absolute error (order of convergence) for the Bernoulli problem}
\label{tab:1}       % Give a unique label
\begin{tabular}{llllllll}
\hline\noalign{\smallskip}
 $k$ &  \centering RK4    &\centering RK6   &\centering DC6RK2/4  &\centering DC6    & BDF6 \\[.5ex]
 4.00e-3 & 1.33e+11        & $--$            &0.54818            & 7.72e-04           & 2.42e-02  \\[.5ex]
 2.00e-3 & 1.1592          & 9.88e+09        &0.2473~(1.14)      & 2.81e-04~(1.45)    & 1.90e-02~(0.35)  \\[.5ex]
 1.00e-3 & 0.353983~(1.71) &    $--$           &4.40e-02~(2.48)    & 4.43e-05~(1.27)    & 1.39e-02~(0.44)  \\[.5ex] 
 1.0e-4  & 1.26e-03~(2.44) & 9.02e-05~       &3.64e-4~(2.08)     & 4.50e-06~(1.27)    & 1.79e-03~(0.89  )\\[.5ex]
 1.0e-5  & 4.91e-08~(4.41) & 2.79e-10~(5.50) &1.16e-09~(5.49)    & 3.92e-09~(3.04)    & 3.15e-06~(2.76) \\[.5ex]
 5.0e-6  & 2.53e-09~(4.27) & 3.51e-12~(6.31) & 9.20e-12~(6.98)   & 1.42e-10~(3.70)    & 1.35e-07~(5.11) \\[.5ex]
 2.5e-6  & 1.42e-10~(4.15) & 1.04e-13~(5.07) & 6.13e-13~(3.90)   & 3.35e-12~(5.54)    & 3.84e-09~(5.13) \\[.5ex]
\hline\noalign{\smallskip}
%Order & \centering 2.01 & \centering 3.99 & 6.03 & 7.94 & 9.90\\[1.2ex]
\end{tabular}
\end{table}

\subsection{Oscillatory problem (see \cite{koyaguerebo2022arbitrary,hull1972comparing, karouma2015class})}
\begin{equation}
\label{66}
u'=\lambda u\cos (t)~,~~~u(0)=1, ~~T=10^6, \lambda=10.
\end{equation}
The exact solution is $u(t)=e^{\lambda \sin(t)}$. The original problem is set with $\lambda =1$ in \cite{hull1972comparing}. Authors in \cite{karouma2015class} solved this problem with Runge-Kutta methods of orders 4 and 8, for $\lambda =2$ and $T=2580\pi$, to `` illustrate the need of higher order methods when a long-term integration problem is considered'', and the current version is due to authors in \cite{koyaguerebo2022arbitrary}. Table \ref{tab:2} gives the absolute error and the order of convergence for each pair of consecutive time steps.

The Magnitude of the exact solution to the oscillatory problem is large, resulting in a large absolute error. While the implicit method $DC6$ is more accurate for this problem at any time step $k$, $DC6RK2/4$ is stable and converges above its proper order. The explicit method $RK6$ requires smaller time steps to produce accurate approximations, and it converges above its proper order, even though it remains less accurate than $DC6RK2/4$ on about two digits. Time-stepping methods $RK4$ and $BDF6$ require sufficiently small time steps to provide accurate approximation to the oscillatory problem.
\begin{table}[!ht]
\caption{Absolute error (order of convergence) for the oscillatory problem}
\label{tab:2}       % Give a unique label
\begin{tabular}{llllllll}
\hline\noalign{\smallskip}
 $k$ &  \centering $RK4$   &\centering$RK6$ &\centering $DC6RK2/4$ &\centering $DC6$ & $BDF6$\\[.75ex]
 5.000e-2 &--            & 2.662e+10     & 9850.859        & 3.2350          & 22026.46\\[.75ex]
 2.500e-2 &3.1e+13      & 2480.0048     & 62.90625~(7.29) & 0.5959~(6.2)    & 22026.46 \\[.75ex]
 1.250e-2 &20354.5     & 18.262~(7.08) & 0.489762~(7.00) & 9.17e-3~(6.0)   & 7659.21~(1.52)\\[.75ex]
 6.250e-3 &454.348~(5.48)&0.1422~(7.00)  & 3.424e-3~(7.16) & 1.4e-4~(5.99)   & 71.00093~(6.75) \\[.75ex]
 3.125e-3 &14.0769~(5.01)&1.12e-3~(6.98) & 3.803e-5~(6.49) & 4.72e-6~($k=\frac{1}{640}$)& 0.551257~(7.00) \\[.75ex]
\hline\noalign{\smallskip}
%Order & \centering 2.01 & \centering 3.99 & 6.03 & 7.94 & 9.90\\[1.2ex]
\end{tabular}
\end{table}

\subsection{Problem B5 (see \cite{koyaguerebo2022arbitrary,enright1975comparing}), stiff with complex eigenvalues of negative real parts}
\begin{equation}
\label{a70}
y'=\begin{bmatrix}-10 & ~~\alpha & ~~0 & ~~0 & 0 & 0\\
               -\alpha & -10     & ~~0 & ~~0 & 0 & 0\\
                   ~~ 0 & ~~0      & -4&~~ 0 & 0 & 0\\
                   ~~ 0 & ~~0      & ~~0 & -1& 0 & 0\\
                   ~~ 0 & ~~0      & ~~0 &~ ~0 & -0.5 & 0\\
                    ~~0 & ~~0      & ~~0& ~~~0~& 0 &-0.1\end{bmatrix}y,~ y(0)=\begin{bmatrix}
                    1\\1\\1\\1\\1\\1
                    \end{bmatrix},~ \alpha=5000,~ T=20.
\end{equation}
This problem, originally set with $\alpha =100$, is an illustration of ODEs resulting from a semidiscretization by finite element methods of parabolic equations (see \cite{stewart1990avoiding}). The current value of $\alpha$ is fixed in \cite{koyaguerebo2022arbitrary} to make the problem more difficult. Table \ref{tab:5} gives the errors for the first component of the approximate solutions, which are similar for the second component. The errors for the other components quickly achieve machine precision. Maximal absolute errors on the six components of the solution for $BDF1$ and $BDF2$ when $k=2\times 10^{-5}$ are, respectively, 1.1857543 and 0.7432313. 

From Table \ref{tab:5} we observe that DC6RK2/4 is better on this problem than the other methods. It achieves its proper order of accuracy and produces very accurate approximations for a large range of the time step $k$ in which $RK6$, $BDF6$ and $DC6$ are unstable. $RK6$, in turn, is more accurate than $BDF6$ and $DC6$.
\begin{table}[!ht]
\caption{Error (order of convergence) for the first component of the solution for $B5$}
\label{tab:5}       % Give a unique label
\begin{tabular}{llllllllll}
\hline\noalign{\smallskip}
 $k$      &\centering $RK4$ & \centering $RK6$ &\centering $DC6RK2/4$ &\centering $DC6$&$BDF6$\\[.75ex]
 4.000e-4 & 1.312598        & --               & 0.9847               & --             &   --  \\[.75ex]
 2.000e-4 & 0.865767~(0.60) & 0.1985           & 8.09e-03~(6.926)     & --             &   --\\[.75ex]
 4.000e-5 & 3.46e-03~(3.435)& 1.101e-5~(6.088) & 5.22e-07~(5.995)     & --             & 2.34e-03\\[.75ex]
 2.000e-5 & 2.16e-04~(3.99)& 1.72e-07~(5.999) & 8.16e-09~(5.998)     & 2.22e-2        & 3.70e-05~(5.98)\\[.75ex]
 5.000e-6 & 8.46e-07~(3.99)& 4.19e-11~(6.001) & 2.04e-12~(5.996)     & 5.59e-6~(6)    & 9.06e-09~(5.997)\\[.75ex]
\noalign{\smallskip}\hline
\end{tabular}
\end{table}

\subsection{Problem E5 (see \cite{enright1975comparing}), stiff with complex eigenvalues of predominantly negative real parts}
\begin{equation}
\label{a71}
 \begin{aligned}
y'_1&=-7.89\times 10^{-10}y_1-1.1\times 10^7y_1y_2&\\
y'_2&=7.89\times 10^{-10}y_1-1.13\times 10^9y_2y_3&\\
y'_3&=7.89\times 10^{-10}y_1-1.1\times 10^7 y_1y_2+1.13\times 10^3y_4-1.13\times 10^9y_2y_3&\\
y'_4&=1.1\times 10^7y_1y_2+1.13\times 10^3y_4&\\
\displaystyle &y(0)=(1.76\times 10^{-3},0;0;0)^t, T=1000.
\end{aligned}
\end{equation}
A reference solution is computed with the solver \texttt{stiff}, taking $rtol=10^8\times atol=10^{-15}$. The solution to this problem has a small magnitude in the region $[1.618\times 10^{-3},1.76\times 10^{-3}]\times [0,1.46\times 10^{-10}]\times [0,8.27\times 10^{-12}]\times [0,1.38\times 10^{-10}]$ of the real plane and the eigenvalues of the Jacobian matrix $dF(y)$ along the solution curve belong to the region  
$[-20490,3.68\times 10^{-12}]\times [-9.17\times 10^{-5},9.17\times 10^{-5}]$ of the complex plane. Table \ref{tab:4} gives the absolute errors and order of accuracy for the four components of the approximate solutions. For the stepsize $k=100$, which corresponds to a uniform partition of the solution interval $[0,1000]$ by eleven points, the maximal absolute errors on the four components of the solution are $8.32\times 10^{-9}$ and $5.45\times 10^{-6}$, respectively, for $DC6$ and $BDF6$. For $k=10$, the maximal absolute errors on the four components of the solution are $5.26\times 10^{-13}$ and $3.52\times 10^{-11}$, respectively, for the latter two methods. We notice that $DC6$ and $BDF6$ achieve machine accuracy in a larger range of stepsize $k$ where the explicit methods are unstable. As a consequence, we restrict the table to absolute errors related to the explicit methods (see \cite{koyaguerebo2022arbitrary} for detailed results about implicit $DC$ methods related to this problem).

From Table \ref{tab:4}, we deduce that implicit methods $DC6$ and $BDF6$ are way more efficient for the problem $E5$ than explicit methods $DC6RK2/4$ and $RK6$ since, for example, $DC6RK2/4$ requires stepsizes about 40000 times smaller to achieve the level of satisfactory accuracy corresponding to the stepzise $k=10$ by $BDF6$ and $BDF6$. Nevertheless, timesteps required by $DC6RK2/4$ to produce accurate approximation are not extremely small. The order of accuracy of the latter method is not observed since the method get in the region of asymptotic convergence with machine accuracy. $RK4$ and $RK6$ are unstable in a large range of the stepsize $k$ for which $DC6RK2/4$ is accurate. 

\begin{table}[!ht]
\caption{ Absolute error (order of convergence) for the problem E5}
\label{tab:4}
\begin{tabular}{lllllllll}
\hline\noalign{\smallskip}~ $k$ &  \centering 1/4000\quad & \centering\quad 1/7000& \centering\quad 1/7500& \quad 1/8000  \\ 
\hline\noalign{\smallskip}
\multirow{4}{6em}{$DC6RK2/4$ }
&7.23e-15 &1.48e-14 &2.03e-14 & 2.27e-14~\\
&6.69e-22 &2.09e-21 &2.49e-21 & 2.64e-21~\\
&1.00e-22 &2.84e-23 &2.38e-23 & 1.60e-23~\\
&6.69e-22 &8.47e-22 &1.26e-21 & 1.65e-21~\\
%\hline\noalign{\smallskip}  
\hline\noalign{\smallskip}   
\multirow{4}{4em}{$RK6$}
& \qquad --&\qquad --& 3.07e-16 &8.28e-16~(-7.40)  \\
& \qquad-- & \qquad--& 1.80e-22 &1.50e-22~(1.34) \\
& \qquad-- & \qquad--& 6.21e-20 &1.31e-23~(63.37) \\
& \qquad --& \qquad--& 6.21e-20 &2.72e-22~(40.66) \\
\hline\noalign{\smallskip}   
\multirow{4}{4em}{$RK4$}
& \qquad --&\qquad --& 2.92e-16 &5.52e-16~(-9.8)  \\
& \qquad-- & \qquad--& 6.41e-21 &6.38e-22~(0.07) \\
& \qquad-- & \qquad--& 1.51e-19 &4.10e-23(127.3) \\
& \qquad --& \qquad--& 1.51e-19 &6.35e-22~(84.8) \\
\noalign{\smallskip}\hline
\end{tabular}
\end{table}

\subsection{ Robertson (1966) \cite{wanner1996solving}, stiff with real negative eigenvalues}
\begin{equation}
\label{69} \begin{aligned}
y'_1&=-0.04y_1+10^{4}y_2y_3\\
y'_2&=0.04y_1-10^{4}y_2y_3-3.10^{7}y_2^2\\
y'_3&=3.10^{7}y_2^2\\
y(0)&=(1,0,0)^t,~~T=10^5.
\end{aligned}
\end{equation}
This is one of the three problems considered stiffest in \cite{wanner1996solving}. Numerical results from the reference \cite{koyaguerebo2022arbitrary} show that the solution belongs to the region $[1.78\times 10^{-2},1.00]\times [0,3.58\times 10^{-5}]\times [0,0.983]$ of the real plane, and the eigenvalues of the Jacobian $dF(y)$ along the solution curve belong to the real interval $[-9825.744,0]$. We compute a reference solution with the solver \texttt{stiff} from Scilab, taking $rtol=100\times atol=10^{-15}$. Table \ref{tab:6} gives absolute errors and orders of accuracy for each component of the solution. For $k=1/2400$, the approximate solution at time $t\simeq 1339.844$, computed with the implemented $RK6$, is $-2.14\times 10^{+33}, -6.43\times 10^{+36}, 6.43\times 10^{+36}$, respectively, for the three components of the solution.

From Table \ref{tab:6}, we deduce that the implicit methods $DC6$ and $BDF6$ are way better efficient on Robertson's problem than the explicit methods, even though approximate solutions from $BDF6$ stagnate around the accurate one obtained from the stepsize $k=1/2$. Nevertheless, $DC6RK2/4$ gives accurate approximate solutions for time steps that are not too small. The order of accuracy is not observed due essentially to the time steps chosen outside its region of asymptotic convergence $k<2/9825.744$. The implemented $RK4$ and $RK6$ face difficulties on the Robertson problem.
\begin{table}[!h]
\caption{Absolute error (order of convergence) for Robertson problem }
\label{tab:6}     
\renewcommand{\arraystretch}{1}\begin{tabular}{llllllll}
%\hline & \multicolumn{6}{|c|}{Global error}\\
\hline\noalign{\smallskip}~ $k$ &\centering RK4 &\centering RK6&\centering DC6RK2/4&\centering DC6 &BDF6\\[1.2ex] 
\hline\noalign{\smallskip}
\multirow{3}{3.5em}{0.5}
   &  --    &   -- &  --       & 2.08e-6 & 4.11e-6 \\
   &  --    &   -- &  --       & 2.08e-6 & 5.33e-9 \\
   &  --    &   -- &  --       & 1.02e-7 & 4.11e-6 \\
\hline\noalign{\smallskip}
\multirow{3}{3.5em}{1/600} 
& --   & --  &  --            & 1.0e-12~(8.6)  & 5.57e-09~(1.15) \\
& --   & --  &  --            & 2e-16~(27.7)   & 8.84e-13~(1.52) \\
& --   & --  &  --            & 1e-12~(15.3)   & 2.64e-08~(0.88)  \\
\hline\noalign{\smallskip}
\multirow{3}{3.5em}{1/1200}
&-- & -- &  --                & 8.33e-13~(0.3) &1.38e-09~(2.01) \\
&-- & -- &  --                & 0.      ~(1.3) &8.39e-15~(6.75) \\
&-- &--  &  --                & 4.01e-13~(1.3) &5.29e-08~(-1.00)\\
\hline\noalign{\smallskip}
\multirow{3}{3.5em}{1/1800}
&--       & -- &1.20e-11       &  &2.07e-09 \\
&--       & -- &2.60e-17       &  &3.18e-14 \\
&--       &--  &7.08e-10               &  &7.93e-08 \\
\hline\noalign{\smallskip}
\multirow{3}{3.5em}{1/2400}
&--       & --  &1.67e-11      &  &2.77e-09 \\
&--       & --  &5.59e-17      &  &1.52e-14 \\
&--       & --  &9.71e-10      &  &1.06e-07 \\
%\hline\noalign{\smallskip}
\noalign{\smallskip}\hline
\end{tabular}
\end{table} 

\subsection{ van der Pol oscillator (see \cite{enright1975comparing,shampine1981evaluation}), stiff, arbitrary complex eigenvalues}

\begin{equation}
\label{a72}
 \begin{aligned}
\displaystyle y'_1&=y_2\\
\displaystyle y'_2&=\mu (1-y_1^2)y_2-y_1\\
\displaystyle y_1(0)&=2,~~y_2(0)=0, T=3000,\mu=1000.
\end{aligned}
\end{equation}
This problem was initially proposed with $T=1$ and $\mu=5$ in \cite{enright1975comparing}. The actual version results from a suggestion by Shampine \cite{shampine1981evaluation}. Numerical results from the reference \cite{koyaguerebo2022arbitrary} show that the solution belongs to the region $[-2, 2.000073]\times [-1323.04,1231.35]$ of the real plane, and the eigenvalues along the solution curve belong to the region $[-3000.29,1123.17]\times [-1158.48,1158.48]$ of the complex plane.  We compute a reference solution with the solver \texttt{stiff} from Scilab, taking $rtol=10atol=10^{-15}$. Table \ref{tab:8} gives the absolute errors and orders of accuracy. For $k=3.75\times 10^{-5}$, the absolute errors on the two components of the approximate solutions are 3.242 and 34.56 for $BDF1$ and 3.24 and 28.54 for $BDF2$. 

From Table \ref{tab:8}, we observe that, as for the problem B5, DC6RK2/4 is better on this problem than the other methods. It achieves its proper order of accuracy and produces very accurate approximations for a large range of the time step $k$ in which approximate solution by $BDF6$ and $DC6$ are not accurate. $RK6$, in turn, is more accurate than $BDF6$ and $DC6$.

\begin{table}[!h]
\caption{Absolute error for the van der Pol's equation }
\label{tab:8}     
\renewcommand{\arraystretch}{1}\begin{tabular}{llllllll}
%\hline & \multicolumn{6}{|c|}{Global error}\\
\hline\noalign{\smallskip}~ $k$&\centering RK4&\centering RK6 &\centering DC6RK2/4&\centering DC6 & BDF6 \\[.2ex] 
\hline\noalign{\smallskip}
\multirow{2}{3.65em}{3.00e-4}
&3.012945    & 2.1269     &0.1229  &3.0938  &3.0265  \\[.2ex]
&1310.151   & 673.95     &269.54  &673.95  &673.93\\[.2ex] 
\hline\noalign{\smallskip}
\multirow{2}{3.65em}{1.50e-4}
&2.958087~(0.02) &2.55e-1(3.06) &8.42e-3~(3.86) &3.0443 &2.9909~(0.017) \\[.2ex]
&1141.658~(0.19) &400.61~(0.75)   &16.8618~(3.99) &673.94 &1296.3~(-0.943)\\[.2ex] 
\hline\noalign{\smallskip}
\multirow{2}{3.65em}{7.50e-5}
&0.684591~(2.11)& 5.10e-3~(5.64) &1.71e-4~(5.62)  &3.02865&6.55e-2 ~(5.00)\\[.2ex]
&1215.283~(-0.1)& 10.102~(5.31)  &3.39e-1~(5.63)  &673.934&123.021 ~(3.39) \\[.2ex] 

\hline\noalign{\smallskip}
\multirow{2}{3.65em}{3.75e-5}
&8.99e-02~(2.92) &9.32e-5~(5.88) &1.70e-6~(6.65) &2.9440~(0.02)  &1.90e-3~(5.10)\\[.2ex]
&89.63520~(3.76) &1.85e-1~(5.76) &6.84e-4~(8.95) &1320.6~(-0.29) &7.58e-1~(7.34) \\[.2ex] 
\noalign{\smallskip}\hline
\end{tabular}
\end{table}

\section{Application to reaction-diffusion equations}
\label{sec:7}
As an application of the hybrid deferred correction method to PDEs, we consider three one dimensional test problems related to reaction-diffusion equations, which are deduced from \cite{koyaguerebo2023arbitrary,owolabi2014higher}. The reaction-diffusion equations, taking the form

\begin{equation}
\label{cyr1} \left\lbrace 
\begin{aligned}
 u_t-\mu u_{xx}+f(x,t,u)&=0, \quad (x,t)\in [x_0,x_f]\times \left( 0,T\right),\\
u(x,0)&=u_0(x), \quad x\in [x_0,x_f]
\end{aligned}
\right. 
\end{equation}
are supposed to be enforced with homogeneous Dirichlet or Neumann boundary conditions (this can be obtained from inhomogeneous cases by changing the unknown function $u$ by $\mathbf{u}=u-\varphi$, where, for example, $\varphi(x,t)=(1-x/xf)u(x_0,t)+(x/x_f)u(x_f,t)$ for DBC and $\varphi(x,t)=(x-0.5x^2/xf)u_x(x_0,t)+(0.5x^2/x_f)u_x(x_f,t)$ for NBC when $0=x_0<x_f$). We discretize (\ref{cyr1}) in space on the uniform grid points $x_0<x_1<\cdots<x_M=x_f$, $x_j=x_0+jh$ and $h=(x_f-x_0)/M$, for an integer $M$, using sixth order finite difference approximations that take into account the boundary conditions. It results the Cauchy problem
 
 \begin{equation}
\label{cyr2} \left\lbrace 
\begin{aligned}
 U_h'+\frac{\mu}{180 h^2} \mathcal{M} U_h+F(U_h,t)&=0, \quad t\in \left( 0,T\right),\\
U_h(0) \mbox{ given in terms of $u_0$}.
\end{aligned}
\right. 
\end{equation}
In the case of Neumann boundary conditions (NBC), $U_h(t)$ and $F(U_h(t),t)$ are column vectors where the $j-$th components are, respectively, $u_j(t)$ and $f(x_j,t,u_j(t))$, $j=0,1,\cdots,M$, and $\mathcal{M}$ is a matrix of size $(M+1)\times (M+1)$ that is equal to the following matrix, denoted $\mathcal{A}$:
$$\left( 
\begin{array}{rrrrrrrrrrrrrr}
360&\;\;-\frac{9958}{7} &6077 &-15126 &\;21290 &-18310 &9609 &\;-2842 &\frac{2552}{7} &0& \cdots&0&\\
-126&70 &486&-855&670&-324&90&-11&0&0&\cdots&0\\
11&-214&378&-130&-85&54&-16&2&0&0&\cdots&0\\
-2&27&-270&490&-270&27&-2&0&0&0&\cdots&0\\
0&-2&27&-270&490&-270&27&-2&0&0&\cdots&0\\
0&0&-2&27&-270&490&-270&27&-2&0&\cdots&0\\
\vdots&&\ddots&\ddots&\ddots& \ddots& \ddots& \ddots&\ddots&\ddots&\ddots&\vdots\\
0&\cdots & 0&-2&27&-270&490&-270&27&-2&0&0\\
0&\cdots &0&0&-2&27&-270&490&-270&27&-2&0\\
0&\cdots &0&0&0&-2&27&-270&490&-270&27&-2\\
0&\cdots&0 &0&2&-16&54&-85&-130&378&-214&11\\
0&\cdots&0 &0&-11&90&-324&670&-855&486&70&-126\\
0&\cdots&0 &\frac{2552}{7}&-2842&9609&-18310&21290 &-15126&6077&-\frac{9958}{7}&360
\end{array}\right).
$$
For Dirichlet boundary conditions (DBC), $U_h(t)$ and $F(U_h(t),t)$ are column vectors where $j-$th components are, respectively, $u_j(t)$ and $f(x_j,t,u_j(t))$, $j=1,2,\cdots,M-1$, and $\mathcal{M}$ is a matrix of size $(M-1)\times (M-1)$, denoted $\mathcal{B}$:
$$\mathcal{B}=\begin{pmatrix}
70&486&-855&670&-324&90&-11&0&0& \cdots&0&\\
-214&378&-130&-85&54&-16&2&0&0&\cdots&0\\
27&-270&490&-270&27&-2&0&0&0&\cdots&0\\
-2&27&-270&490&-270&27&-2&0&0&\cdots&0\\
0&-2&27&-270&490&-270&27&-2&0&\cdots&0\\
\vdots&\ddots&\ddots&\ddots&\ddots&\ddots&\ddots&\ddots&\ddots&\ddots&\vdots\\
0&\cdots &0&-2&27&-270&490&-270&27&-2&0\\
0&\cdots &0&0&-2&27&-270&490&-270&27&-2\\
0&\cdots &0&0&0&-2&27&-270&490&-270&27\\
0&\cdots &0&0&2&-16&54&-85&-130&378&-214\\
0&\cdots &0&0&-11&90&-324&670&-855&486&70\\
\end{pmatrix}. $$
The Cauchy problem (\ref{cyr2}) combines larger system and stiffness. The first test problem related to reaction-diffusion equations, which is the Fisher equation, has exact solution and will be used for a quick evaluation of the order of accuracy in space for the discretization (\ref{cyr2}) in its two cases related to Dirichlet and Neumann boundary conditions.

To evaluate the errors, we use the norm

$$\|\left\lbrace U^{n}_h\right\rbrace_n-u_h\| =\max _{0\leq n\leq N}|U^{n}_h-u_h(t_n)|,$$
where $\vert \,\cdot \,\vert $ denotes the Euclidean norm, extracting the solutions essentially at 100 discrete times evenly spread over $[0,T]$. As in the numerical experiments (see Section \ref{sec:6}), we comparatively use RK4, RK6, DC6RK2/4, DC6 and BDF6 to approximate solutions of the reaction-diffusion equations.

\subsection{Fisher equation (see \cite{arora2020meshfree,koyaguerebo2023arbitrary})}
\begin{figure}[h!]
        \begin{minipage}[t]{8cm}
        \centering
        \includegraphics[width=7cm]{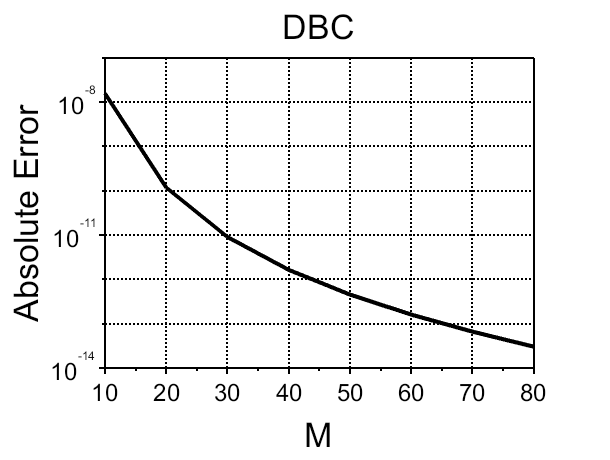}
%        \caption{Absolute error in the space variable for the Fisher equation with DBC, fixing $k=4.00\times 10^{-4.}$}
    \end{minipage}   
            \begin{minipage}[t]{8cm}
        \centering
        \includegraphics[width=7cm]{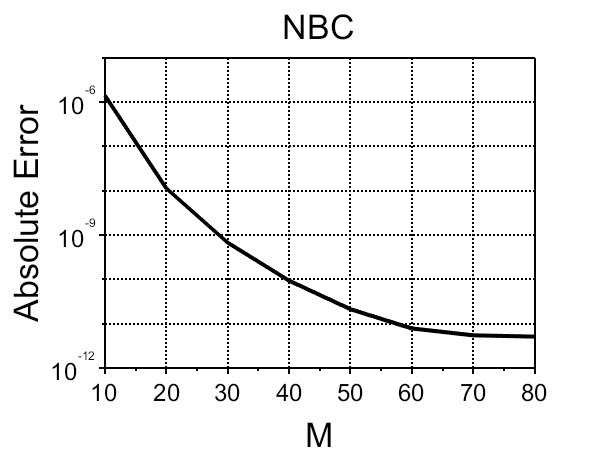}
%        \caption{Absolute error in the space variable for the Fisher equation with NBC, fixing $k=4.00\times 10^{-4.}$}
    \end{minipage}
\label{fig:fisher}       % Give a unique label
\caption{Absolute errors in the space variable for the Fisher equation with DBC (left) and NBC (right), fixing $k=10^{-3}$.}
\end{figure}

\begin{equation}\displaystyle 
\label{f3} \begin{aligned}
u_t-u_{xx}-6u(1-u)&=0,\quad (x,t)\in \Omega\times [0,T],\\
u(x,0)&=(1+e^{x})^{-2}, \quad x\in \Omega.
\end{aligned}
\end{equation}
The exact solution is $u(x,t)=(1+e^{x-5t})^{-2}$ (see \cite{ablowitz1979explicit}). We study the cases of Dirichlet boundary condition (DBC) and Neumann boundary conditions (NBC), knowing that the boundary conditions are given by the exact solution. This problem is used as a test problem in \cite{arora2020meshfree} for $\Omega=[0,1]$ and $T=0.002$. Authors in \cite{koyaguerebo2023arbitrary} solve the problem for the same space interval, but for $T=10$. As in \cite{koyaguerebo2023arbitrary}, we choose $\Omega=[0,1]$ and $T=10$. We consider the discretization given by (\ref{cyr2}), fixing $M=80$. Table \ref{tab:Fisher} gives absolute errors and some orders of accuracy in time, and Figure \ref{fig:fisher} gives the graph of absolute error in the space variable, computed with DCIM6 where $k$ is fixed to $10^{-3}$.

\begin{center}
\begin{table}[!ht]
\caption{Absolute error (order of convergence) for the Fisher equation with $M=80$ fixed}
\label{tab:Fisher}       % Give a unique label
\begin{tabular}{llllllll}
\hline\noalign{\smallskip}
 &        &  & DBC &  & \\[1.2ex]  
\hline\noalign{\smallskip}
 $N$    &  \centering RK4& \centering RK6 &\centering DC6RK2/6  &\centering DC6& BDF6\\
 100    &  --           & --               & --              & 1.42e-4         & 3.72e-03      \\
 1000   &  --           & --               & --              & 1.17e-09~(5.08) & 1.12e-08~(5.52)\\[.7ex]
 10000  &  --           & --               & --              & 3.02e-14~(4.58) & 3.61e-14~(5.49)\\[.7ex]
 70000  &  --           & --               & 3.03e-14        & 2.99e-14~(0.01) & 7.30e+07\\[.7ex]
 120000 &  --           & --               & 2.74e-14        & 3.07e-14(-0.04) & 1.02e-13\\[.7ex]
 140000 &5.41e-14       & 5.37e-14         & 5.22e-14        & 3.07e-14(-0.00) & 1.16e-13(-1.0)\\[.7ex]
\hline\noalign{\smallskip}
 &        & &  NBC &  & \\[1.2ex]
\hline\noalign{\smallskip}
 $N$    &  \centering RK4& \centering RK6 &\centering DC6RK2/6  &\centering DC6& BDF6\\
 100    &  --           & --               & --              & 1.09e-04        & 1.16e-02\\[.7ex]
 1000   &  --           & --               & --              & 4.40e-10~(5.39) & 2.90e-08~(5.60)\\[.7ex]
 100000 &  --           & --               & 5.54e-12        & 5.10e-12        & 2.00e+08\\[.7ex]
 160000 & --            & 7.19e-12         & 5.03e-12        & 5.52e-12        & 6.26e+08 \\[.7ex]
 200000 & 1.09e-11      & 1.05e-11         & 5.17e-12        & 1.01e-11        & 1.37e-10\\[.7ex]
\hline\noalign{\smallskip}
\end{tabular}
\end{table}
\end{center}
%%NBC BDF[r100 r1e3 r1e4 r1e5 r2e5] 0.0116403   2.900D-08   1.143D-11   2.009D+08   1.377D-10
%
%--> [log(r100./r1e3)/log(10) log(r1e3./r1e4)/log(10) log(r1e4./r1e5)/log(10) log(r1e5./r2e5)/log(2)]
%   5.6035839   3.4043993  -19.24498   60.339242
%DCIM

%--> --> [w100 w1e3 w1e4 w1e5 w2e5]   0.0001091   4.407D-10   5.100D-12   4.520D-12   1.060D-11

%--> [log(w100./w1e3)/log(10) log(w1e3./w1e4)/log(10) log(w1e4./w1e5)/log(10) log(w1e5./w2e5)/log(2)]
%
%   5.3936064   1.9365983   0.0524461  -1.2292217

%RK
%--> [v9e4;v1e5;v16e4;v2e5]
%
% ans = [4x3 double]
%
%   0.          0.          5.549D-12
%   0.          0.          5.249D-12
%   0.          7.192D-12   5.039D-12
%   1.095D-11   1.052D-11   5.173D-12

%DBC
%DCRK--> [v7e4; v12e4; v14e4]
%   0.          0.          3.036D-14
%   0.          0.          2.746D-14
%   5.414D-14   5.370D-14   5.224D-14

%-->DCIM2-6 [w100;w1e3;w1e4;w7e4]
%
% ans = [4x3 double]
%
%   0.0070114   0.0005816   0.0001422
%   0.000066    5.619D-08   1.171D-09
%   0.0000007   5.684D-12   3.023D-14
%   1.347D-08   3.137D-14   2.999D-14
%
%--> [log(w100./w1e3)/log(10);log(w1e3./w1e4)/log(10);log(w1e4./w7e4)/log(7/4)]
%
% ans = [3x3 double]
%
%   2.0260257   4.0150036   5.0841605
%   2.0002429   3.9949479   4.5881596
%   6.9544585   9.2915819   0.0139765
   
From Figure \ref{fig:fisher}, we observe a convergence in the space variable with the expected order 6 of accuracy both for the Dirichlet and Neumann boundary conditions. Otherwise, from Table \ref{tab:Fisher}, we observe that the time steps require for accurate approximations of the Fisher equation are not excessively small as it is expected for explicit methods. While implicit methods are better on the Fisher equation with DBC for larger number of steps $N>4400$, the approximations with $DC6RK2/4$ for this problem achieve machine accuracy for $k=1/7000$ and $k=10^{-4}$, respectively, for the DBC and the NBC. In both cases, DC6RK2/4 achieves machine accuracy on a large range of the time step $k$ for which RK4 and RK6 are unstable.

\subsection{Bistable reaction-diffusion equation (see \cite{koyaguerebo2023arbitrary})}
\begin{figure}[!h]%% placement specifier
%% Use \includegraphics command to insert graphic files. Place graphics files in 
%% working directory.
\centering%% For centre alignment of image.
\includegraphics[height=3in, width=4.5in]{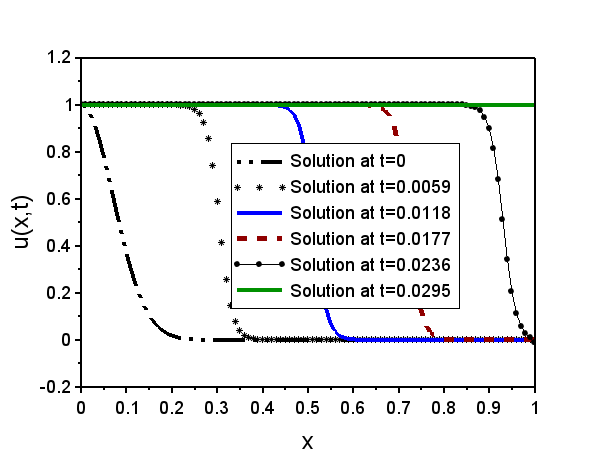}
%% Use \caption command for figure caption and label.
\caption{Reference solution of the bistable reaction-diffusion equation at times $t=0,\; 0.0059,\;0.0118,\;0.0177,\; 0.0236,\; 0.0295$.}
\label{fig2}
%% https://en.wikibooks.org/wiki/LaTeX/Importing_Graphics#Importing_external_graphics
\end{figure}

\begin{equation}\displaystyle 
\label{f1} \begin{aligned}
u_t-u_{xx}+10^4u(u-1)(u-0.25)&=0 \mbox{ in } \Omega\times (0,T),\\
\frac{\partial u}{\partial n}&=0\mbox{ on } \partial \Omega\times (0,T),\\
u(0)=u_0(x)&=e^{-100x^2} \mbox{ in } \Omega.\\
\end{aligned}
\end{equation}
As in \cite{koyaguerebo2023arbitrary}, we choose $\Omega =(0,1)$ and $T=0.0295$. We apply the formula (\ref{cyr2}) for the space discretize, fixing $M=100$. We compute a reference solution of the Cauchy problem with the ODE solver \texttt{stiff} from Scilab, taking $rtol=10^9\times atol=10^{-15}$. The solution is a travelling wave that leaves the domain at about $t=0.0295$ (see Figure \ref{fig2}). The maximal absolute errors in time for $BDF1$ and $BDF2$ when $N=100$ are, respectively, 15.61 and 1.71. For $N=6300$, the erros are 0.279 and $1.12 \times 10^{-3}$, respectively, for $BDF1$ and $BDF2$.

    \begin{table}[!ht]
\caption{Absolute error (order of convergence) for the bistable reaction-diffusion equation}
\label{tab:9}       % Give a unique label
\begin{tabular}{llllllll}
\hline\noalign{\smallskip}
 $N$  &  \centering RK4 & \centering RK6 & \centering DC6RK2/6  &\centering DC6& BDF6\\
 100  &  --           & --               & --              & 2.68e-02       & 0.35197      \\
 500  &  --           & --               & 8.59e-07        & 4.47e-06(5.40) & 40.891(-2.9)\\[.7ex]
 800  &  --           & --               & 2.96e-08(7.16)  & 2.45e-07(6.17) & 30.467(0.62)\\[.7ex]
 1200 &9.28e-06       & 1.84e-06         & 2.05e-09(6.58)  & 2.09e-08(6.07) & 1.45e-06(41.5)\\[.7ex]
 3600 &9.85e-08(4.13) & 1.40e-09(6.53)   & 1.41e-11(4.53)  & 3.45e-11(5.82) & 2.73e-09(5.71)\\[.7ex]
 6300 &1.02e-08(4.04) & 4.55e-11(6.12)   & 1.39e-11(0.02)  & 1.41e-11(1.60) & 1.01e-10(5.88)\\[.7ex]
 \hline\noalign{\smallskip}
%Order & \centering 2.01 & \centering 3.99 & 6.03 & 7.94 & 9.90\\[1.2ex]
\end{tabular}
\end{table}

From Table \ref{tab:9} we observe that, except for $N<500$, $DC6RK2/4$ is more accurate than all the other time-stepping methods, and it produces very accurate approximations for a large range of the time step $k$ in which $RK4$ and $RK6$ are unstable. All the methods achieve their proper order on this problem.

\subsection{The Robertson equation with one dimensional diffusion (see \cite{khaliq2009smoothing, owolabi2014higher} and references therein)}
\begin{equation}\displaystyle 
\label{f2} \begin{aligned}
u_t-\alpha u_{xx}+\tau_1 u-\tau_2 v\,w &=0 \mbox{ in } \Omega\times (0,T),\\
v_t-\beta v_{xx}-\tau_3 u+\tau_4 v\,w+\tau_5 v^2 &=0 \mbox{ in } \Omega\times (0,T),\\
w_t-\gamma w_{xx}-\tau_6 v^2 &=0 \mbox{ in } \Omega\times (0,T),\\
u(0,t)=u(1,t)=1;v(0,t)=v(1,t)=0; w(0,t)=w(1,t)&=0, \quad t\in (0,T),\\
u(x,0)=1+\sin(2\pi x);v(x,0)=0;w(x,0)&=0,\quad x \in \Omega.
\end{aligned}
\end{equation}
As in \cite{owolabi2014higher}, we choose $\Omega =(0,1)$, $\tau_1=\tau_3=4\times 10^{-2}$, $\tau_2=\tau_4=10^{4}$,$\tau_5=\tau_6=3\times 10^{7}$, $\alpha=\beta=\gamma=1$, but we evaluate the solution up to $T=1$. The case where a diffusion coefficient is equal to 10 (see \cite{khaliq2009smoothing}) simply requires a time step $k_0/10$ for the explicit methods, if $k_0$ is the initial time step for the explicit method to produce an accurate approximation in the present case. The discretization of this problem, using the finite differences in (\ref{cyr2}), writes
\begin{equation}
\label{f2b}
\left\lbrace 
\begin{aligned}
U'_h+\frac{1}{180h^2}\mathcal{B}U_h+F_1(U_h,V_h,W_h)&=0, \quad t\in (0,T]\\
V'_h+\frac{1}{180h^2}\mathcal{B}V_h+F_2(U_h,V_h,W_h)&=0, \quad t\in (0,T]\\
W'_h+\frac{1}{180h^2}\mathcal{B}W_h+F_3(U_h,V_h,W_h)&=0, \quad t\in (0,T]\\
U_h(0)&= \begin{bmatrix}
u_0(x_1)-1 &u_0(x_2)-1&\cdots &u_0(x_{M-1}-1)
\end{bmatrix}^T,\\
V_h(0)&= \begin{bmatrix}
v_0(x_1) &v_0(x_2)&\cdots &v_0(x_{M-1})
\end{bmatrix}^T,\\
W_h(0)&= \begin{bmatrix}
w_0(x_1) &w_0(x_2)&\cdots &w_0(x_{M-1})
\end{bmatrix}^T,
\end{aligned}\right. 
\end{equation}

\begin{figure}[h!]
        \begin{minipage}[t]{5cm}
        \centering
        \includegraphics[width=5cm]{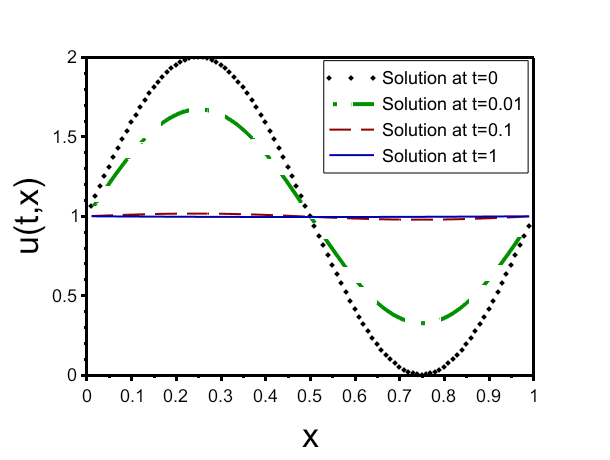}
        %\caption{}
    \end{minipage}   
            \begin{minipage}[t]{5cm}
        \centering
        \includegraphics[width=5cm]{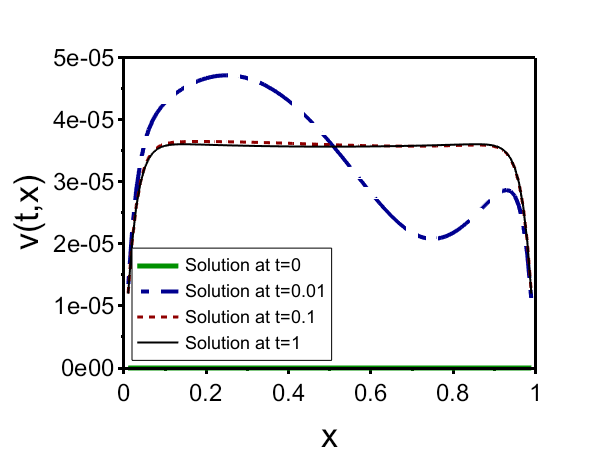}
        %\caption{}
    \end{minipage}
    \begin{minipage}[t]{5cm}
        \centering
        \includegraphics[width=5cm]{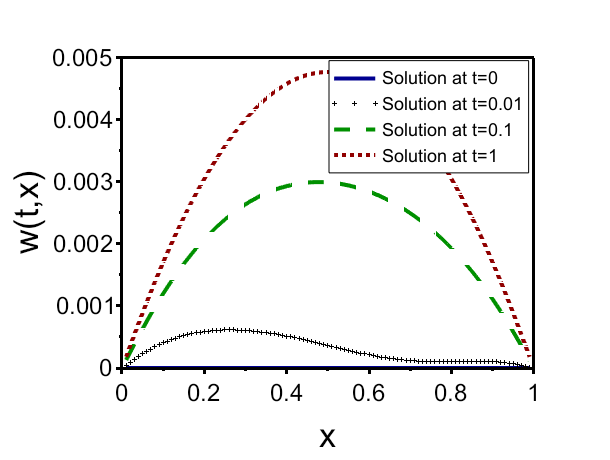}
        %\caption{}
    \end{minipage}   
\caption{The three components of the reference solution for the three-specie system at times $t=0,0.001, 0.01, 0.1,1$.}
\label{fig:Kolade}       % Give a unique label
\end{figure}

We fix $M=100$ and compute a reference solution of the Cauchy problem (\ref{f2b}) with the ODE solver \texttt{stiff} from Scilab, taking $rtol=atol=10^{-16}$. Figure \ref{fig:Kolade} gives a graph of the three components of the reference solution at times $t=0,0.001,0.01,0.1,1$, and Table \ref{tab:10} gives absolute errors and orders of accuracy for each component of the solution.
% For $N=100$, the maximal absolute error on the three components of the approximate solutions by BDF1 and BDF2 are, respectively, 0.250 and 0.102. For $N=12000$, the absolute errors on the three components of the approximate solutions are $3.48\times 10^{-3}$ and $3.03\times 10^{-5}$, respectively, for BDF1 and BDF2.  The Runge-Kutta methods $RK4$ and $RK6$ fail for solving this problem when $N\leq 21000$.

From Table \ref{tab:10}, we deduce that the implicit methods perform better on the system (\ref{f2b}) than explicit methods, even though approximate solutions by $BDF6$ are unstable around the time step $k=1/11000$. Nevertheless, $DC6RK2/4$ gives accurate approximate solutions for time steps that are not excessively small and provides the more accurate solution in its convergence region.

\begin{table}[!h]
\caption{Absolute error (order of convergence) for the three-specie system}
\label{tab:10}     
\renewcommand{\arraystretch}{1}\begin{tabular}{llllllll}
%\hline & \multicolumn{6}{|c|}{Global error}\\
\hline\noalign{\smallskip}~ $N$ &\centering RK4&\centering RK6&\centering DC6RK2/4 & \centering DC6 & BDF6\\[1.2ex] 
\hline\noalign{\smallskip}
\multirow{3}{3.5em}{100}
&--            &  --          &     --         &4.58e-05         &1.03e-03 \\
&--            &  --          &     --         &4.58e-05         &1.03e-03  \\
&--            &  --          &     --         &3.11e-05         &1.33e-05  \\
\hline\noalign{\smallskip}
\multirow{3}{3.5em}{1000}
&--            &  --          &     --         &4.54e-08~(3.00)  &7.47e-07~(3.13) \\
&--            &  --          &     --         &4.54e-08~(3.00)  &7.47e-07~(3.13)  \\
&--            &  --          &     --         &1.21e-08~(3.40)  &7.97e-06~(0.22)  \\
\hline\noalign{\smallskip}
\multirow{3}{3.5em}{11000}
&--            &  --          & 4.89e-03       &1.08e-13~(5.39)  &-- \\
&--            &  --          &4.89e-03        &1.08e-13~(5.39)  &--\\
&--            &  --          &1.02e-02        &1.14e-14~(5.78)  &--\\
\hline\noalign{\smallskip}

\multirow{3}{3.5em}{12000}
&--            &  --          &1.60e-13    &7.17e-13 &-- \\
&--            &  --          &1.60e-13    &7.16e-13  &-- \\
&--            &  --          &1.60e-14    &7.16e-13  &--  \\
\hline\noalign{\smallskip}
\multirow{3}{3.5em}{24000}
&2.51e-13  &  1.65e-13   &1.71e-13)  &7.19e-13   &4.09e-11 \\
&2.51e-13  &  1.65e-13   &1.71e-13   &7.16e-13   &4.09e-11 \\
&2.89e-14  &  1.23e-14   &1.23e-14   &7.17e-13   &3.70e-11  \\
\noalign{\smallskip}\hline
\end{tabular}
\end{table} 

\section{General Observations}
\label{sec:8}
The hybrid deferred correction method involving explicit midpoint rule the classical fourth-order Runge-Kutta method, $DC6RK2/4$, is stable on each of the test problems and does not require extremely small time steps for accurate numerical solutions. The method achieves its proper order on almost each of the test problems and handles problems addressing nonlinearity and long-term integration. Moreover, as shown on the problem B5 and the van der Pol oscillator,  $DC6RK2/4$ is better than the implicit DC methods and the $BDF$ family on problems of type (\ref{a1}) where Jacobian matrices $d_uF(t,u(t))$ along the solution curve have eigenvalues $\lambda=\lambda_1+i\lambda_2$ where $-\lambda_1 \leq \vert \lambda_2\vert$, due to its region of absolute stability that contains a significant part of the imaginary axis (for stability) together with its high order and smaller error constant (for accuracy). We refer to \cite[Section 6.7]{koyaguerebo2022arbitrary} and \cite{stewart1990avoiding} for inefficiency of BDF methods and the DC methods for the implicit midpoint rule on such stiff problems.

Each step of the scheme $DC6RK2/4$ requires five sub-steps of $RK4$ in the correction steps, leading to a computational cost of about 5.25 times more than that of $RK4$ to compute approximate solutions on the same grid points. Nevertheless, the higher order of $DC6RK2/4$ together with its larger region of absolute stability allows this method to be essentially more efficient than $RK4$. In fact, the largest region of absolute stability of $DC6RK2/4$ with respect to that of $RK4$ and $RK6$ allows this method to achieve machine accuracy (case of the Robertson equation and some reaction-diffusion equations) or to produce very accurate approximate solutions in a large range of the time step $k$ for which $RK4$ and $RK6$ are unstable. On the other hand, the effect of the higher order of $DC6RK2/4$ with respect to $RK4$ can be observed, for example, on the problem B5 for which absolute errors for $DC6RK2/4$ and $RK4$ are, respectively, $5.22\times 10^{-7}$ and $8.46\times 10^{-7}$ while the number of function evaluations for these two methods are, respectively, 1050000 and 1600000. Similarly, on the bistable reaction-diffusion equation, when $N=1200$ for $DC6RK2/4$ and $N=6300$ for $RK4$, theses two methods use the same number of function evaluations while the absolute error for $DC6RK2/4$ is $2.05\times 10^{-9}$ and that of $RK4$ is $1.05\times 10^{-8}$. Otherwise, contrarily to the IDC in \cite{IntegralDC2010}, that requires at least 33, 22 or 27 function evaluations per step for sixth-order IDC schemes constructed on $RK2$, $RK3$ and $RK4$, respectively, $DC6RK2/4$ requires only 21 function evaluations per step and is very simple to implement. Moreover, $DC6RK2/4$ is embedded since each approximate solution at a node $t_{n+1}$ by $DC6RK2/4$ with a step $k$ uses an approximate solution at $t_{n+1}$ by $RK4$ in the correction step, computed with the uniform step $k/5$ on the interval $[t_n,t_{n+1}]$.

\end{document}